\documentclass[12pt]{article}

\usepackage{graphicx}
\usepackage{amsmath, amsthm}
\usepackage{amsfonts}
\usepackage{amssymb}
\usepackage{url}
\usepackage{verbatim}
\usepackage{color}
\definecolor{red}{rgb}{1,0,0}

\definecolor{blu}{rgb}{0,0,1}

\usepackage[mathlines]{lineno}
\usepackage{mathrsfs}
\definecolor{qqqqff}{rgb}{0.,0.,1.}

\def\noi{\noindent}

\def\wh#1{\widehat{#1}}
\newlength{\sdgwidth}
\def\sdiv#1{\settowidth{\sdgwidth}{$#1$}
\mathop{\ooalign{$\mkern1.6mu\overline{\protect\raisebox{0.4pt}{$
\phantom{\mkern-1.6mu#1\kern-0.5\sdgwidth\mkern-.8mu}$}}$\hidewidth\protect\cr$
\kern0.5\sdgwidth\mkern0.8mu\mathring{\protect\raisebox{-.25pt}{$
\phantom{\kern-0.5\sdgwidth#1\kern-0.5\sdgwidth}$}}
\mkern-0.8mu\kern0.5\sdgwidth$\protect\cr$#1$\protect\cr$
\hidewidth\overline{\protect\raisebox{0.4pt}{$
\phantom{\kern-0.5\sdgwidth\mkern-2.2mu#1}$}}$}}\nolimits}

\def\pdk#1{\gamma_{P,#1}}

\newcommand{\kpd}{\gamma_{P,k}}
\newcommand{\kf}{Z_k}
\newcommand{\sk}{{\mathscr P}_{G,k}}
\newcommand{\hk}{{\mathscr P}_{H,k}}
\newcommand{\fk}{{\mathscr F}_{G,k}}

\newcommand{\Z}{\operatorname{Z}}

\newtheorem{thm}{Theorem}[section]
\newtheorem{cor}[thm]{Corollary}
\newtheorem{lem}[thm]{Lemma}
\newtheorem{prop}[thm]{Proposition}

\newtheorem{obs}[thm]{Observation}

\theoremstyle{definition}

\theoremstyle{definition}
\newtheorem{defn}[thm]{Definition}
\theoremstyle{definition}
\newtheorem{ex}[thm]{Example}

\newcommand{\bit}{\begin{itemize}}
\newcommand{\eit}{\end{itemize}}
\newcommand{\ben}{\begin{enumerate}}
\newcommand{\een}{\end{enumerate}}
\newcommand{\beq}{\begin{equation}}
\newcommand{\eeq}{\end{equation}}
\newcommand{\bea}{\begin{eqnarray*}}
\newcommand{\eea}{\end{eqnarray*}}
\newcommand{\bpf}{\begin{proof}}
\newcommand{\epf}{\end{proof}\ms}
\newcommand{\ms}{\medskip}

\newcommand{\lc}{\left\lceil}
\newcommand{\rc}{\right\rceil}
\newcommand{\lf}{\left\lfloor}
\newcommand{\rf}{\right\rfloor}

\title{The relationship between $k$-forcing and $k$-power domination}

\author{ Daniela Ferrero\thanks{Department of Mathematics, Texas State University, San Marcos, TX 78666, USA (dferrero@txstate.edu)}
\and Leslie Hogben\thanks{Department of Mathematics, Iowa State University, Ames, IA 50011, USA (hogben@iastate.edu) and American Institute of Mathematics, 600 E. Brokaw Road, San Jose, CA 95112, USA (hogben@aimath.org).}\and Franklin H.J. Kenter\thanks{Department of Mathematics, United States Naval Academy, Annapolis, MD 21402, USA (franklin.kenter@gmail.com).}\and Michael Young \thanks{Department of Mathematics, Iowa State University, Ames, IA 50011, USA (myoung@iastate.edu)}}

\begin{document}

\maketitle

\begin{abstract} 
Zero forcing and power domination are iterative processes on graphs where an initial set of vertices are observed, and additional vertices become observed based on some rules. In both cases, the goal is to eventually observe the entire graph using the fewest number of initial vertices. Chang et al. introduced $k$-power domination in [Generalized power domination in graphs, {\it Discrete Applied Math.} 160 (2012) 1691-1698] as a generalization of power domination and standard graph domination.  Independently, Amos et al. defined $k$-forcing in [Upper bounds on the $k$-forcing number of a graph, {\it Discrete Applied Math.} 181 (2015) 1-10] to generalize zero forcing. In this paper, we combine the study of $k$-forcing and $k$-power domination, providing a new approach to analyze both processes. We give a relationship between the $k$-forcing and the $k$-power domination numbers of a graph that bounds one in terms of the other.  We also obtain results using the contraction of subgraphs that allow the parallel computation of $k$-forcing and $k$-power dominating sets.
\end{abstract}

\noi {\bf Keywords} $k$-power domination, $k$-forcing, subgraph contraction, Sierpi\'nski graphs

\noi{\bf AMS subject classification} 05C69, 05C50 

\section{Introduction}

Zero forcing was introduced as a process to obtain an upper bound for the maximum nullity of real symmetric matrices whose nonzero pattern of off-diagonal entries is described by a given graph  \cite{AIM08}. The minimum rank problem was motivated by the inverse eigenvalue problem of a graph. Independently, zero forcing was introduced by mathematical physicists studying quantum systems \cite{graphinfect}. Since its introduction, zero forcing has attracted the attention of a large number of researchers who find the concept useful to model processes in a broad range of disciplines. The need for a uniform framework for the analysis of the diverse processes where the notion of zero forcing appears led to the introduction of a generalization of zero forcing called $k$-forcing  \cite{ACDP14}.  

Amos et al. proposed $k$-forcing in \cite{ACDP14} as the following graph coloring game. Assume the vertices of a graph are colored in two colors, say white and blue. Iteratively apply the following color change rule: if $u$ is a blue vertex with at most $k$ white neighbors, then change the color of all the neighbors of $u$ to blue. Once this rule does not change the color of any vertex, if all vertices are blue, the original set of blue vertices is a $k$-forcing set of $G$. The original zero forcing is $1$-forcing under this definition. Because the problem of deciding whether a graph admits a $1$-forcing set of a given maximum size is NP-complete even if restricted to planar graphs \cite[Theorem 2.3.1]{A10}, the general problem of finding forcing sets cannot be solved algorithmically for large graphs without the development of further theoretical tools.

Power domination was introduced by Haynes et al. in \cite{HHHH02} when using graph models to study the monitoring process of electrical power networks. When a power network is modeled by a graph, a power dominating set provides the locations where monitoring devices (Phase Measurement Units, or PMUs for short) can be placed in order to monitor the power network. Finding optimal PMU placements is an important practical problem  in electrical engineering due to the cost of PMUs and network size.  Although power domination is substantially different from standard graph domination, the notion of $k$-power domination was proposed as a generalization of both power domination ($k=1$) and standard graph domination ($k=0$) \cite{CDMR12}. 

Chang et al. defined $k$-power domination in \cite{CDMR12} using sets of {\it observed} vertices.  Given a graph $G$ and a set of vertices $S$, initially all vertices in $S$ and their neighbors are observed; all other vertices are unobserved. Iteratively apply the following propagation rule: if there exists an observed vertex $u$ that has $k$ or fewer unobserved neighbors, then all the neighbors of $u$ are observed. Once this rule does not  produce any additional observed vertices, if all vertices of $G$ are observed, $S$ is a $k$-power dominating set of $G$.  Many problems outside graph theory can be formulated in terms of minimum $k$-power dominating sets \cite{CDMR12} so methods to obtain them are highly desired. An algorithmic approach has been attempted, but the problem of deciding if a graphs admits a $k$-power dominating set of a given maximum size is NP-complete \cite{CDMR12}. 

Although  $k$-forcing  and $k$-power domination have been studied independently, an in-depth analysis of $k$-power domination leads to the study of $k$-forcing. Indeed, after the initial step in which a set observes itself and its neighbors, the observation process in $k$-power domination proceeds exactly as the color changing process in $k$-forcing. The aim of this paper is to establish a precise connection between $k$-forcing and $k$-power domination to facilitate the transference of results, proofs, and methods	 between them, and ultimately to advance research on both problems.

Throughout this paper we work on $k$-forcing and $k$-power domination concurrently, using results in one process as stepping stones for results in the other one. In Section 2 we present the definitions and notation that we use in the rest of the paper. In Section 3 we give some core results and remarks that we use in the sections that follow.

In Section 4 we examine the effect of subgraph contraction in $k$-power domination and $k$-forcing. We obtain upper and lower bounds for the change in the $k$-power domination number produced by the contraction of a subgraph. Note that the contraction of a subgraph can increase or decrease its $k$-power domination number. In particular, we prove that the contraction of subgraphs of small degree can  change the $k$-power domination number by at most one. In this section we also propose a way to decompose a graph in order to bound its $k$-power domination number in terms of that of smaller subgraphs. This can allow computation of $k$-power dominating sets to run in parallel. We also give the analogous results for $k$-forcing.

In Section 5 we present a lower bound for the  $k$-power domination number of a graph in terms of its $k$-forcing number.  This bound generalizes a known result for $k=1$ that gives the only lower bound for the power domination number of an arbitrary graph available so far \cite{REUF2015}. As an application, we find an upper bound for the $k$-forcing number of a graph in terms of its maximum degree.\vspace{-3pt}


\section{Definitions and notation}\vspace{-3pt}

A  {\em graph} is an ordered pair $G=(V,E)$ where $V=V(G)$ is a finite nonempty set of {\em  vertices}  and $E=E(G)$ is a set of unordered pairs
of distinct vertices  called {\em edges}  (i.e., in this work graphs are simple and undirected). The {\em order} of $G$ is $|G|:=|V(G)|$. Two vertices $u$ and $v$ are {\em adjacent} or {\em neighbors} in $G$ if $\{u,v\} \in E(G)$. The {\em (open) neighborhood} of a vertex $v$ is the set $N_G(v) = \{u \in V: \{u,v\} \in E\}$, and the {\em closed neighborhood} of $v$ is the set $N_G[v] = N_G(v) \cup \{ v \}$. Similarly, for any set of vertices $S$, $N_G(S) = \cup_{v \in S} N_G(v)$ and $N_G[S] = \cup_{v \in  S} N_G[v]$. The {\em degree} of a vertex $v$ is $\deg _G (v):=|N(v)|$. The {\em maximum} and {\em minimum degree} of $G$ are $\Delta(G) := \max \{\deg _G (v):v\in V\}$ and $\delta(G) := \min \{\deg _G (v):v\in V\}$, respectively; a graph $G$ is {\em regular} if $\delta(G)=\Delta(G)$. We will omit the subscript $G$ when the graph $G$ is clear from the context.

A path joining $u,v\in V$ is a sequence of vertices $u=x_0,x_1,\ldots ,x_r=v$ such that $\{x_i,x_{i+1}\}\in E$ for each $i=0,\ldots , r-1$. A graph $G$ is {\it connected} if there is a path joining every pair of different vertices. If a graph is not connected, each maximal connected subgraph is a {\it component} of $G$. In this paper, $c(G)$ denotes the {\it number of components} of $G$ and $G_1,\ldots ,G_{c(G)}$ denote the components of $G$. Most of the results in this work are given for connected graphs, since if a graph is not connected, we can apply the results to each component.

If $X$ is a set of vertices of $G$,  the  {\it subgraph induced by $X$  (in $G$)} is denoted as $G[X]$; it has vertex set $X$ and edge set $\{\{u,v\}\in E:u,v\in X\}$. The graph $G-X$ is defined as $G[V\setminus X]$. The {\it contraction of $X$ in $G$} is the graph $G/ X$ obtained by adding a vertex $v_X$ to $G-X$ with  $N_{G/X}(v_X)=N_G[X]\setminus X$. Note that $G/X$ does not require $G[X]$ to be connected whereas the standard use of graph contraction does.

In a graph $G=(V,E)$, consider an arbitrary coloring of its vertices in two colors, say blue and white, and let $T$ denote the set of blue vertices. The color changing process in $k$-forcing can be formally described by associating to $T$ the family of sets $(\fk ^i (T))_{i\geq 0}$ recursively defined by the following rules.\vspace{-3pt}

\begin{itemize}
\item[1.] $\fk ^0(T)=T$,\vspace{-3pt}
\item[2.] $\fk ^{i+1} (T)=\fk ^i (T)\cup \{N(v): v \in  \fk ^i (T) \mbox{ and } 1\le | N_G(v)\setminus \fk ^i (T)|\le k \}$, for $i\geq 0$.
\end{itemize}

A set $T\subseteq V$ is a {\em $k$-forcing set} of  $G$ if there is an integer $t$ such that
$\fk ^t(T)=V$. A {\em minimum $k$-forcing set} is a $k$-forcing set of minimum cardinality.  The {\em $k$-forcing number}  of $G$ is 
the cardinality of a minimum $k$-forcing set and is denoted by $\kf (G)$. If  $v\in \fk ^{i} (T)$ and $| N(v)\setminus \fk ^i (T)|\le k$ then $v$ is said to $k${\it-force} (or simply {\it force} if $k$ is clear from the context) every vertex in $N(v)\setminus \fk ^i (T)$.

Let $k$ be a nonnegative integer. The definition of $k$-power domination on a graph $G$ will be given in terms of a family of sets, $(\sk ^i (S))_{i\geq 0}$,  associated to each set of vertices $S$ in $G$.\vspace{-3pt}

\begin{itemize}
\item[1.] $\sk ^0(S)=N[S]$,
\item[2.] $\sk ^{i+1} (S)= \sk ^i (S)\cup \{N(v):  v \in  \sk ^i (S) \mbox{ and } 1\le | N_G(v)\setminus \sk ^i (S)|\le k \}$, for $i \ge 0$.
\end{itemize}

A set $S\subseteq V$ is a {\em $k$-power dominating set} of  $G$ if there is an integer $\ell$ such that
$\sk ^{\ell}(S)=V$. A {\em minimum $k$-power dominating set} is a $k$-power dominating set of minimum cardinality. The {\em $k$-power domination number} of $G$ is the cardinality of a minimum $k$-power dominating set and is denoted by $\kpd (G)$.  

Next we recall the definition of standard graph domination. A vertex $v$  \emph{dominates} all vertices in $N_G[v]$.  A set $S\subseteq V$ is a \emph{dominating set} of $G$ if $N_G[S]=V$.  The minimum cardinality of a dominating set is the \emph{domination number} of $G$, denoted by $\gamma (G)$.

Note that $1$-forcing coincides with zero forcing \cite{ACDP14}, while $1$-power domination is exactly power domination and $0$-power domination coincides with domination \cite{CDMR12}.

\section{Preliminaries}

The following observations follow directly from the definitions of $k$-power domination and $k$-forcing, and provide the initial connection between both concepts.

\begin{obs}\label{definitionsets}
In any graph $G$, if $T$ is a $k$-forcing set, all sets $( \fk ^i (T))_{i\geq 0}$ are $k$-forcing sets of $G$; if $S$ is a $k$-power dominating of $G$, the sets $(\sk ^i (S))_{i\geq 0}$ are also $k$-forcing sets of $G$.
\end{obs}

\begin{obs}\label{basic}
In any  graph $G$, if $T$ is a $k$-forcing set of $G$ then $T$ is also a $k$-power dominating set. The converse is not necessarily true, but $S$ is a $k$-power dominating set if and only if $N[S]$ is a $k$-forcing set. As a consequence, $\gamma _{P,k}(G)\leq Z_k(G)\leq \gamma _{P,k}(G) (\Delta (G)+1)$.
\end{obs}

\begin{obs}\label{neighbors}
In a graph $G$, 
$S\subsetneq V(G)$ is a $k$-power dominating set of  $G$ if and only if $N[S]\setminus S$ is a $k$-forcing set of $G-S$.
\end{obs}

Note that given a graph $G=(V,E)$ and $S\subseteq X \subseteq V$, it is possible that for some $x\in X$, $\deg_{G[X]}(x)<\deg_G(x)$. Therefore, the $k$-power domination process starting with $S$ in $G$ is different from the one starting with $S$ in $G[X]$. As a consequence,  $S$ being a $k$-power dominating set of $G[X]$ does not imply that $S$ can $k$-observe all vertices in $X$ when propagating in $G$. Analogously, if $T\subseteq X \subseteq V$ then $T$ being a $k$-forcing set of $G[X]$ does not imply that $T$ can $k$-force $X$ in $G$. This observation motivates the following definitions.

\begin{defn}\label{defkfsubset}
Let $G=(V,E)$ be a graph and let $A\subseteq X \subseteq V$. We say that $A$ is a $k${\em -forcing set  of $X$ in $G$} if there exists a nonnegative integer $t$ such that $X\subseteq { \mathscr F}^t_{G,k}(A)$.
\end{defn}

\begin{defn}\label{defpdsubset}
Let $G=(V,E)$ be a graph and let $A\subseteq X \subseteq V$. We say that $A$ is a $k${\em -power dominating set  of $X$ in $G$} if there exists a nonnegative integer $\ell$ such that $X\subseteq { \mathscr P}^{\ell}_{G,k}(A)$.
\end{defn}

The proofs of the next results are straightforward, and are omitted.

\begin{lem}\label{transitive}
Let $T$ be a $k$-forcing set of a graph $G$. Let $A\subseteq T$.
\begin{itemize}
\item [1)] If $A$ is $k$-forcing set of $T$ in $G$, then $A$ is a $k$-forcing set of $G$;
\item [2)] If $A$ is $k$-power dominating set of $T$ in $G$, then $A$ is a $k$-power dominating set of $G$.
\end{itemize}
\end{lem}

\begin{lem}\label{transitive2}
Let $S$ be a $k$-power dominating set of a graph $G$. 
Let $A\subseteq S$.
\begin{itemize}
\item [1)] If $A$ is $k$-forcing set of $N[S]$ in $G$, then $A$ is a $k$-forcing set of $G$;
\item [2)] If $A$ is $k$-power dominating set of $N[S]$ in $G$, then $A$ is a $k$-power dominating set of $G$.
\end{itemize}
\end{lem}


\begin{lem}\label{onevertex}
Let $G=(V,E)$ be a graph and $X\subseteq V$ such that $G[X]$ is connected and $\deg_G(x)\leq k+1$ for every $x\in X$. Let $u$ be an arbitrary vertex in $X$.  Then  $\{u\}$ is a (minimum) $k$-power dominating set of $N[X]$ in $G$. In addition,  if $\deg_G(u)\leq k$, then $\{u\}$ is also a (minimum) $k$-forcing set of $N[X]$ in $G$.
\end{lem}

\bpf 
If $S=\{u\}$, then $\sk ^0(S)=N[u]$ and $\fk ^0(S)=\{u\}$. Since $\deg_G(x)\leq k+1$ for every $x\in X$,  $x\ne u$ has at most $k$ unobserved neighbors when $x$ is observed.  Thus,  $\sk ^i(S)=N[\sk ^{i-1}(S)]$  for every integer $i\geq 1$. Since $G[X]$ is connected, there exists an integer $r\geq 1$ such that $X\subseteq \sk ^r(S)$. Once all vertices in $X$ are observed, each of them can have at most $k$ unobserved neighbors, so such a vertex can observe any unobserved neighbors. Thus,  $N[X]\subseteq \sk ^{r+1}(S)$ and $S$ is a $k$-power dominating set of $N[X]$ in $G$. Now suppose
$\deg_G(u)\leq k$. Then  $\fk ^1(S)=N[S]$ and the argument proceeds as before. 
\epf

The following result follows immediately from Lemma \ref{onevertex}, but is already known for $k$-power domination \cite[Lemma 7]{CDMR12};  a slightly weaker version   for $k$-forcing is given in \cite[Proposition 2.3]{ACDP14}.

\begin{cor}\label{onevertexnumbers}
Let $G$ be a connected graph. If $\Delta (G) \leq k+1$, then $\pdk k (G)=1$; if in addition $\delta (G)\leq k$,  then $Z_k(G)=1$. 
\end{cor}

When $G[X]$ is not connected, we apply Lemma \ref{onevertex} in each of its components and obtain the following result.

\begin{cor}\label{onevertexbordernc}
Let $G=(V,E)$ be a connected graph, $X\subseteq V$ and $u_j\in V(G[X]_j)$ for every $j=1,\ldots ,c(G[X])$. Let $S=\{u_1, \ldots ,u_{c(G[X])} \}$. If $\deg_G(x)\leq k+1$ for every $x\in X$,  then  $S$ is a minimum $k$-power dominating set of  $N[X]$ in $G$;  if in addition $\deg_G(u_j)\leq k$ for every $j=1,\ldots ,c(G[X])$,  then $S$ is a minimum $k$-forcing set of $N[X]$ in $G$.
\end{cor}

\bpf By Lemma \ref{onevertex}, for every $j=1,\ldots , c(G[X])$,  $\{u_j\}$ is a $k$-power dominating set of $G[X]_j$ in $G$.  Thus, $S$ is a $k$-power dominating set of $N[X]$ in $G$. Since every $k$-power dominating set of $N[X]$ must have at least one vertex in each component of $G[X]$ and $|S|= c(G[X])$ we conclude that $S$ is a minimum $k$-power dominating set of $N[X]$ in $G$. The argument for $k$-forcing is analogous.
\epf

\begin{lem}\label{4}
Let $G=(V_G,E_G)$ and $H=(V_H,E_H)$ be two graphs. 
Let $A\subseteq V_G$ and  $B\subseteq V_H$ such that (i) $G-A=H-B$ and (ii) $N_G[A]\setminus A=N_H[B]\setminus B$. Then
\begin{itemize}
\item [1)] $A$ is a $k$-power dominating set of $G$ if and only if $B$ is a $k$-power dominating set of $H$;
\item [2)] $N_G[A]$ is a $k$-forcing set of $G$ if and only if $N_H[B]$ is a $k$-forcing set of $H$.
\end{itemize}
\end{lem}

\bpf
1) If $A$ is a $k$-power dominating set of $G$, then $N_G[A]\setminus A$ is a $k$-forcing set of $G-A$ by Observation \ref{neighbors}. Since $G-A=H-B$ and $N_G[A]\setminus A=N_H[B]\setminus B$,  we substitute $N_G[A]\setminus A$ and $G-A$ with $N_H[B]\setminus B$  and $H-B$, respectively, and obtain that $N_H[B]\setminus B$ is a $k$-forcing set of $H-B$. By Observation \ref{neighbors},  $B$ is a $k$-power dominating set of $H$. 

2) If $N_G[A]$ is a $k$-forcing set of $G$, then $A$ is a $k$-power dominating set of $G$ by Observation \ref{basic}. Using 1) we conclude that $B$ is a $k$-power dominating set of $H$, and by Observation \ref{basic} we conclude that $N_H[B]$ is a $k$-forcing set of $H$.
\epf

\begin{cor}\label{5}
Let $G=(V_G,E_G)$ and $H=(V_H,E_H)$ be two graphs. 
Let $A\subseteq V_G$ and $B\subseteq V_H$ such that (i) $G-A=H-B$ and (ii) $N_G[A]\setminus A=N_H[B]\setminus B$.  Let $P\subseteq V_G\setminus A=V_H \setminus B$. Then
\begin{itemize}
\item [1)] $A\cup P$ is a $k$-power dominating set of $G$ if and only if $B\cup P$ is a $k$-power dominating set of $H$;
\item [2)] $N_G[A]\cup P$ is a $k$-forcing set of $G$ if and only if $N_G[B]\cup P$ is a $k$-forcing set of $H$.
\end{itemize}
\end{cor}

\bpf
Define $A'=A\cup P$ and $B'=B\cup P$. Then $G-A'=H-B'$ and $N_G[A']\setminus A'=N_H[B']\setminus B'$, so we can apply Lemma \ref{4} with $G$, $H$, $A'$ and $B'$.
\epf

While all the previous results include analogous statements for $k$-forcing and a $k$-power domination, the following lemma does not have a $k$-forcing analog. 

\begin{lem}\label{highdegkPSD}{\rm \cite[Lemma 9]{CDMR12}} 
If $G$ is connected and $\Delta (G)\geq k+2$, then there exists a minimum $k$-power dominating set $S$ such that $\deg (v) \ge k+2$ for all $v \in S$.
\end{lem}

To see that there is no $k$-forcing analog to Lemma \ref{highdegkPSD} it is sufficient to consider $K_{1,n}$, the complete bipartite graph with one vertex in one part and $n$ vertices in the other. As shown in \cite{ACDP14}, if $n>k$ every minimum $k$-forcing set contains at least one vertex of degree $1$. 


\section{Graph contraction}

\begin{defn}\label{1}
Let $G$ be a graph and let $X\subseteq V(G)$. Define $\widehat{X}$ to be the graph obtained from $G[X]$ by attaching to each one of its vertices as many pendent vertices as its number of neighbors in $G-X$. 
\end{defn}

\begin{lem}\label{3}
Let $G$ be a connected graph and let $X\subseteq V(G)$.  
There exists $S\subseteq X$ such that $S$ is a minimum $k$-power dominating set of $\widehat {X}$.
\end{lem}

\bpf  Suppose first that $\Delta (\widehat X)\leq k+1$.  Then by Lemma \ref{onevertex} any one vertex of $X$ is a $k$-power dominating set for $N[X]=\wh X$; a one vertex $k$-power dominating set is necessarily minimum. Now assume $\Delta (\widehat X)\geq k+2$.
By definition of  $\widehat X$, $\deg_{\widehat X}(u)=1$ for every $u\in V(\widehat X)\setminus X$. Since $\Delta (\widehat X)\geq k+2$, by Lemma \ref{highdegkPSD} there exists a minimum $k$-power dominating set  $S$ of $\widehat X$ that contains only vertices in $X$. 
\epf

For the same reasons why there is no $k$-forcing analog to Lemma \ref{highdegkPSD}, there is no $k$-forcing analog to Lemma \ref{3}. Indeed, if $x\in V(G)$ and $\deg_G(x)\geq k+1$, a minimum $k$-forcing set of $\widehat {\{x\}}$ must contain a vertex of degree 1.

\begin{lem}\label{new3}
Let $G$ be a connected graph and let $X\subseteq V(G)$.  If $S\subseteq X$ is a minimum $k$-power dominating set of $\widehat {X}$, then $S$ is a $k$-power dominating set of $N_G[X]$ in $G$.
\end{lem}

\bpf
Each vertex in $V(\widehat X)\setminus X$ arises from a vertex $y\notin X$ that is a neighbor of a vertex $x\in X$. For every $x\in X$, let $N_x$ denote the (possibly empty) set of neighbors of $x$ in $V(\widehat X)\setminus X$ (i.e.,  $N_x=N_{\widehat X}(x)\setminus X$) and let $N'_x=N_G(x)\setminus X$. Since $S\subseteq X$ and $\deg_{\widehat X}(u)=1$ for every $u\in V(\widehat X)\setminus X$, none of the vertices in $N_x$ can be observed before $x$ is observed, and moreover, all vertices in $N_x$ are observed simultaneously.  Since for every $x\in V(\widehat X)$, $\deg_{\widehat X}(x)=\deg _G(x)$, the only difference between the $k$-power domination process starting with $S$ in $\widehat X$ and the one starting with $S$ in $G$ is that when the vertices in $N_x$ are observed in $\widehat X$,  the unobserved vertices in $N'_x$ become observed in $G$.  The reason why some vertices in $N'_x$ could have been observed earlier is that a vertex in $G-X$ could have more than one neighbor in $X$ so $(N'_x)_{x\in X}$ are not necessarily disjoint. Since for every $w\in N_G[X]\setminus X$ there exists $x\in X$ such that $w\in N'_x$, all vertices in $N_G[X]$ are observed.
\epf

\begin{thm}\label{subgraphpd}
Let $G=(V,E)$ be a connected graph. If $X\subseteq V$, \[\gamma _{P,k} (G/X)-1\leq \gamma _{P,k} (G)\leq \gamma _{P,k} (G/X)+ \gamma _{P,k}(\widehat X)\] and both bounds are tight.
\end{thm}

\bpf

Let $H=G/X$.  By Lemma \ref{3} there exists $\widehat P\subseteq X$ such that $\widehat P$ is a minimum $k$-power dominating set of $\widehat X$ and by Lemma \ref{new3}, $\widehat P$ is  also a $k$-power dominating set of $N_G[X]$ in $G$. 

To prove the upper bound we show that if $P$ is a $k$-power dominating set of $H$, then $(P\setminus \{v_X\}) \cup \widehat P$ is a $k$-power dominating set of $G$.\footnote {Note that whether $v_X\notin P$ or $v_X\in P$ does not affect the conclusion, since in any case $|S|\leq |P|+ |\widehat P| =\gamma _{P,k} (H)+ \gamma _{P,k}(\widehat X) $; we only exclude $v_X$ from $S$ to guarantee $S\subseteq V(G).$}  
Since $\widehat P$ is a $k$-power dominating set of $N_G[X]$ in $G$,  clearly $(P\setminus \{v_X\}) \cup \widehat P$ is a $k$-power dominating set of $N_G[P\setminus \{v_X\}]\cup N_G[X]=N_G[P\setminus \{v_X\}\cup X]$ in $G$. We will prove that $(P\setminus \{v_X\})\cup X$ is a $k$-power dominating set of $G$, which by Lemma \ref{transitive2} suffices to conclude that $(P\setminus \{v_X\}) \cup \widehat P$ is a $k$-power dominating set of $G$. Let $A=X$ and $B=\{v_X\}$. Since $H=G/X$, $G- A=H-B$ and $N_G[A]\setminus A=N_H[B]\setminus B$, we apply  Corollary \ref{5} and conclude that $(P\setminus \{v_X\})\cup A$ is a $k$-power dominating set of $G$ if and only if $(P\setminus \{v_X\})\cup B$ is a $k$-power dominating set of $H$. Since $B=\{v_X\}$, $(P\setminus \{v_X\})\cup B= P$ and $P$ is a $k$-power dominating set of $H$,   $(P\setminus \{v_X\})\cup A=(P\setminus \{v_X\})\cup X$ is a $k$-power dominating set of $G$.

To prove the lower bound, we show that if $S$ is a minimum $k$-power dominating set of $G$, then  $(S\setminus X)\cup \{v_X\}$ is a $k$-power dominating set of $H$.  As above, let $A=X$ and $B=\{v_X\}$ so $G- A=H-B$ and $N_G[A]\setminus A=N_H[B]\setminus B$.  Then we apply Corollary \ref{5} to conclude that  $(S\setminus X)\cup A$ is a $k$-power dominating set of $G$ if and only if $(S\setminus X)\cup B$ is a $k$-power dominating set of $H$. Since $X=A$, then $(S\setminus X)\cup A=S$ and it is a $k$-power dominating set of $G$. Then $(S\setminus X)\cup B=(S\setminus X)\cup \{v_X\}$ is a $k$-power dominating set of $H$.  Thus, $\gamma _{P,k} (G/X)\leq |(S\setminus X)\cup \{v_X\}|\leq |S|+1=\gamma _{P,k} (G)+1$.

To prove the upper bound is tight, for each integer $q\geq k$ we define a graph $U_q$ and a set $X\subseteq V(U_q)$ such that $ \gamma _{P,k} (U_q)= \gamma _{P,k} (U_q/X)+ \gamma _{P,k}(\widehat {X})$ (see Figure 1). Consider two disjoint copies of $K_{q+2}$, say $G$ and $G'$,  and vertices $x\in V(G)$ and $y\in V(G')$. Construct $U_q$ by adding the edge $e=\{x,y\}$ and define  $X=V(G')\setminus \{y\}$.  Then    $\gamma _{P,k}( {U_q}) =2$, $\gamma _{P,k}( \widehat {X}) =1$, and  $\gamma _{P,k}(U_q / X) =1$.

\begin{figure}[h!] \begin{center}
\scalebox{.5}{\includegraphics{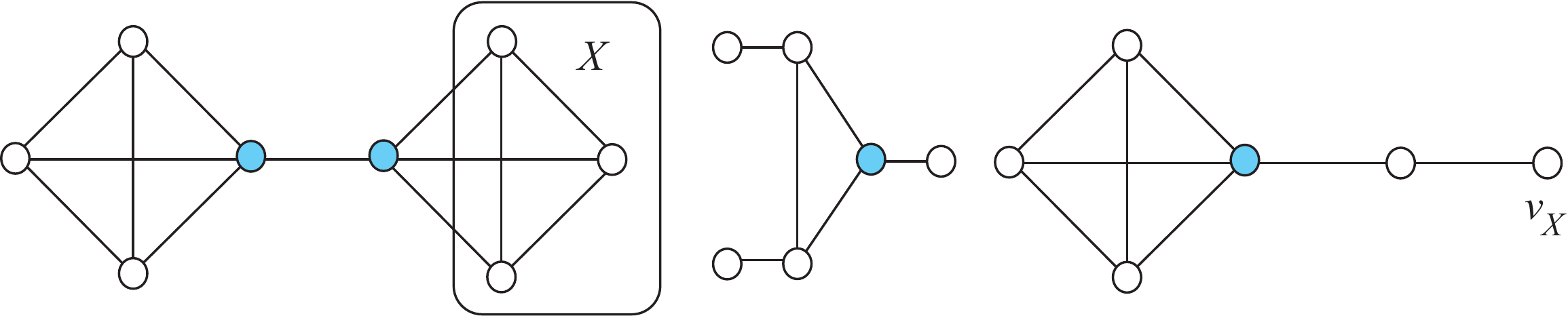}}\\
$U_2$\qquad\qquad\qquad\qquad\qquad $\wh X$ \qquad\qquad\qquad\qquad $U_2/X$
\caption{The graphs $U_2$, 
$\wh X$, and $U_2/X$  defined in Theorem \ref{subgraphpd} are shown.  In each case, a minimum $2$-power dominating set 
is indicated by coloring.\vspace{-10pt}}

\end{center}\end{figure}

\begin{figure}[h!] \begin{center}
\scalebox{.5}{\includegraphics{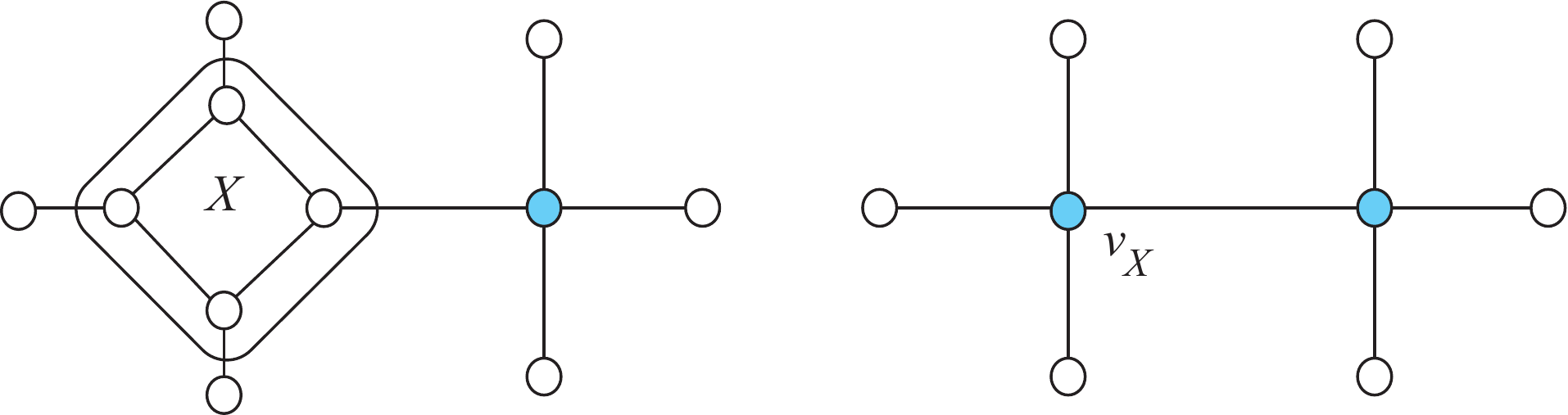}}\\
$L_2$\qquad\qquad\qquad\qquad\qquad  $L_2/X$
\caption{The graphs $L_2$ and $L_2/X$ defined in Theorem \ref{subgraphpd} are shown. 
In each case, a minimum $2$-power dominating set 
is indicated by coloring.\vspace{-10pt}}
\end{center}\end{figure}

To show the lower  bound is tight, for each integer $q\geq k$ we define a graph $L_q$ and a set $X\subseteq V(L_q)$ such that $\gamma _{P,k} (L_q/X)-1= \gamma _{P,k} (L_q)$ (see Figure 2). Assume first that $k\geq 2$.  Construct $L_q$ starting with a cycle of length $2q$ with vertices $v_1,\ldots ,v_{2q}$. Attach a pendent vertex to each vertex $v_i$, for $i=1,\ldots ,2q$. Then attach $q+1$ pendent vertices to the pendent neighbor of $v_1$, so  $\gamma _{P,k}(L_q) =1$.  For  $X=\{v_1,\ldots ,v_{2p}\}$, $\gamma _{P,k}(L_q / X)=2$. Now suppose $k=1$, and begin with a path of length $6$ with vertices $v_0, \ldots ,v_6$. Construct $L_q$ by attaching $q$ pendent vertices to $v_0$. If $X=\{v_1,v_3,v_5\}$, then $\gamma _{P,1}(L_q) =1$ and  $\gamma _{P,1}(L_q / X)=2$.
\epf

The next example shows that it is possible to find a graph $G$ and a subgraph $X$ for which the gap between $\gamma _{P,k}(G) $ and $\gamma _{P,k}(G/X)$ is arbitrarily large.

\begin{ex}  Given a positive integer $c$, define $T_{k,c}$ as the tree obtained by adding $k+2$ leaves to each leaf of $K_{1,c}$. If $X$ is the set of all vertices of degree greater than one in $T_{k,c}$, then $\gamma _{P,k}(T_{k,c})= \gamma _{P,k}(\widehat X)=c$ and  $\gamma _{P,k}(T_{k,c}/X)=1$.

\end{ex}

\begin{cor}\label{subgraphcontraction}
	Let $G=(V,E)$ be a connected graph. Let $X\subseteq V$ such that $\deg_G(x)\leq k+1$ for every $x\in X$. Then \[ \gamma _{P,k}(G/X) -1 \le \gamma _{P,k}(G) \leq  \gamma _{P,k}(G/X) +c(G[X]).\] 
\end{cor}

\bpf
It is sufficient to show that $\gamma _{P,k}(\widehat X)\leq c(G[X])$. Observe that  $\deg_G(x)\leq k+1$ for every $x\in X$ implies that  $\deg_{\widehat X}(x)\leq k+1$ for every $x\in V(\widehat X)$. 
By Corollary \ref{onevertexbordernc}  there exists a $k$-power dominating set of $N_{\widehat X}[X]=V(\widehat X)$ in $G$ with cardinality $c(\widehat X)=c(G[X])$, so $\gamma _{P,k}(\widehat X)\leq c(G[X])$.
\epf

\begin{cor}\label{nice}
	Let $G=(V,E)$ be a connected graph. Let $X\subseteq V$ such that $G[X]$ is connected and $\deg_G(x)\leq k+1$ for every $x\in X$.Then \[ \gamma _{P,k}(G/X) -1 \le \gamma _{P,k}(G) \leq  \gamma _{P,k}(G/X) +1.\] 
\end{cor}

\begin{prop}\label{improvement}
	Let $G=(V,E)$ be a connected graph. Let $X\subseteq V$ such that $G[X]$ is connected and $\deg_G(x)\leq 2$ for every $x\in X$. Then $ \gamma _{P,1}(G/X) \leq  \gamma _{P,1}(G)$ and this bound it tight.
	\end{prop}

\bpf

Since $\Delta (G)\leq 2$ implies that $G$ itself is a path or a cycle, without loss of generality we can assume $\Delta (G)\geq 3$. By Lemma \ref{highdegkPSD}, there exists a minimum $k$-power dominating set $S$ of $G$ such that $\deg_G(u)\geq 3$ for every $u \in S$, so  $S\subseteq V\setminus X$. We prove that $S$ is also a $k$-power dominating set of $H$. 

As in the proof of Theorem \ref{subgraphpd}, $S\cup\{v_X\}$ is a $k$-power dominating set of $H$. Then by Observation \ref{basic}, $N_H[S\cup \{ v_X \}]=N_H[S]\cup N_H[v_X]$ is a $k$-forcing set of $H$.

Note that $S\subseteq V\setminus X$ implies $\sk ^{0} (S)\setminus X=\hk ^{0} (S)\setminus \{v_X\}$ and as long as $\sk ^{i} (S)\subseteq V\setminus X$, $\sk ^{i} (S)=\hk ^{i} (S)$.  Let $x$ be a vertex of $X$ that is observed first (meaning that no vertex of $X$ has been observed earlier), and let $y$ be the vertex in $G-X$ that dominates or forces  $x$ at time $t$ ($x\in \sk ^{0} (S)$ or $x\in \sk ^{t} (S)\setminus \sk ^{t-1} (S)$ for $t\ge 1$).  Since $\deg_G(y)\geq \deg_H(y)$, $y$ can also dominate or force in $H$.  Thus $v_X \in \hk ^{t} (S)$. Since $\deg_H(v_X)\leq 2$, it takes at most one additional application of $k$-forcing to observe all vertices in $N_H[v_X]$, so $N_H[v_X]\subseteq \hk ^{t+1} (S)$.

Since $N_H[S]= \hk ^{0} (S)\subseteq  \hk ^{t+1} (S)$ and  $N_H[v_X]\subseteq \hk ^{t+1} (S)$, then $N_H[S\cup \{ v_X \}]=N_H[S]\cup N_H[v_X]\subseteq  \hk ^{t+1} (S)$. Moreover, since $N_H[S\cup \{ v_X \}]$ is a $k$-forcing set of $H$, so is $\hk ^{t+1} (S)$ and therefore, $S$ is a $k$-power dominating set of $H$.

To prove the tightness, observe that for $n\ge 3$, contracting the set of all vertices of degree $2$ in the path $P_{n}$ of order $n$ produces the path $P_3$. Now, $\gamma _{P,1}(P_{n})=\gamma _{P,1}(P_{3})=1$.
\epf

Due to the computational complexity of the $k$-power domination problem, efficient algorithms to approximate of optimal $k$-power dominating sets are of practical importance. Theorem \ref{subgraphpd} could help in the parallel search for $k$-power dominating sets. The following result provides a theoretical framework to study practical uses of graph decomposition as a tool for the parallel computation of $k$-power dominating sets.

\begin{thm}\label{partitionpd}
Let $G=(V,E)$ be a  connected graph and let $P_1,\ldots ,P_r$ be a partition of $V$. Then \[\gamma _{P,k}(G)\leq \sum_{i=1}^r \gamma _{P,k}(\widehat {P_i}).\]
\end{thm}

\bpf
By Lemma \ref{3}, for every $i=1,\ldots ,r$ there exists $S_i\subseteq P_i$ such that $S_i$ is a minimum $k$-power dominating set of $\widehat {P_i}$. By Lemma \ref{new3}, $S_i$ is also a $k$-power dominating set of $N_G[P_i]$ in $G$, and as a consequence, $S=\cup _{i=1}^r S_i$ is a $k$-power dominating set of $G$. Then $\gamma _{P,k}(G)\leq |S|\leq  \sum_{i=1}^r |S_i| =\sum_{i=1}^r \gamma _{P,k}(\widehat {P_i})$. 
\epf

To prove that the bound in Theorem \ref{partitionpd} is tight we will use the family of Sierpi\'nski graphs whose definition we recall, using the notation in \cite{DK14}. Given two positive integers $n$ and $p$ the {\em Sierpi\'nski graph } $S_p^n$ has as vertices all $n$-tuples of integers in $\{0,1,\ldots ,p-1\}$ denoted as $s_ns_{n-1}\cdots s_2 s_1$. Two vertices $s_n\cdots s_1$ and $t_n\cdots  t_1$ are adjacent in $S_p^n$ if and only if there exists an $r$ with $1\leq r\leq n$ such that

\begin{itemize}
\item[i)] $s_i=t_i$ for every $i\in \{r+1,\ldots ,n\}$,
\item[ii)] $s_r\not=t_r$, and
\item[iii)] $s_i=t_r$  and $t_i=s_r$ for every $i\in \{1,\ldots ,r-1\}$
\end{itemize}

The definition of Sierpi\'nski graphs implies that   $S_p^1=K_p$ and if $n\geq 2$, $S_p^n$ has $p^{n-i}$ induced copies of $S_p^{i}$.  
Moreover, the vertices in each of those copies coincide in the $n-i$ leftmost digits $s_n\cdots s_{i+1}$  \cite{DK14}.  If $s$ is a $(n-i)$-tuple of integers in $\{0,\ldots ,p-1\}$, let $sS_p^i$ denote the set of vertices of $S_p^n$ whose leftmost $n-i$ digits coincide with $s$. For simplicity, we use $S_p^n[s]$ to denote $S_p^n[sS_p^{i}]$ (the subgraph induced by $sS_p^i$ in $S_p^n$), so $S_p^n[s]$ is isomorphic to $S_p^{i}$.

\begin{lem}
Given integers $n\ge 4$, $k\geq 1$ and $p\geq k+2$,  let  $s$ be a $(n-3)$-tuple of integers in $\{0,\ldots ,p-1\}$.  Then $\gamma_{P,k}(\widehat {sS_p^3})= \gamma_{P,k}(S_p^{3})$.
\end{lem}

\bpf
Fix $k\ge 1$, $p\ge k+2$, and $n\ge 4$.  
We begin by determining the degree of a vertex $sxyz$ in $S^n_p$ and in $S_p^n[s]\cong S_p^3$.
The definition of the Sierpi\'nski graph $S_p^n$ implies that vertices of the form $a^n$ have degree $p-1$ and  all the other vertices have degree $p$ in $S_p^n$. If $xyz$ is nonconstant, then vertex $sxyz$ has degree $p$ in both $S_p^n$ and $S_p^n[s]$, so $sxyz$ does not have any pendent vertices added  in $\widehat {sS_p^3}$.  Now consider the constant sequence $aaa$.  If $s\ne a^{n-3}$, then vertex $saaa$ has degree $p$ in  $S_p^n$ but degree $p-1$ in $S_p^n[s]$, so one leaf is added to $saaa$  in $\widehat {sS_p^3}$. If $s= a^{n-3}$, then vertex $saaa=a^n$ has degree $p-1$ in both $S_p^n$ and  $S_p^n[s]$, so $saaa$ is unchanged in $\widehat {sS_p^3}$.

If $s\ne a^{n-3}$ for any $a\in \{0,\ldots ,p-1\}$, then $\widehat {sS_p^3}$ is obtained from $S_p^n[s]$ by attaching one pendent vertex  to each of the $p$ vertices of the form $saaa$ for $a\in \{0,\ldots ,p-1\}$, so that every vertex of $\widehat {sS_p^3}$ has degree $p$; denote this graph by $G_1$.   If $s=a^{n-3}$ for some $a\in \{0,\ldots ,p-1\}$, then $\widehat {sS_p^3}$ is obtained from $S_p^n[s]$ by attaching one pendent vertex to each of the $p-1$ vertices of the form $sbbb$ for $b\ne a$, so that every vertex  of $\widehat {sS_p^3}$ except $a^n$ has degree $p$; denote this graph by $G_2$.
Observe that the only difference between $G_1$ and $G_2$ is that $G_2$ is missing one leaf.

To show that $\gamma_{P,k}(G_{i})= \gamma_{P,k}({S_p^3})$ for $i=1,2$, we first we prove $\gamma_{P,k}(G_{i})\leq \gamma_{P,k}({S_p^3})$ by showing that if $P$ is a  $k$-power dominating set of $S_p^{3}$, then $P$ is also a $k$-power dominating set of $G_{i}$. For each pendent vertex $x$ in $G_{i}$, let $u_x$ denote its only neighbor (in $G_{i}$). Then $u_x$ is a vertex of degree $p-1$ in $S_p^3$ and therefore, in $S_p^3$ it is labeled as $a^3$ for some $a\in \{0,\ldots ,p-1\}$. If $u_x\in P$, then $x\in N_{G_{i}}[P]={{\mathscr P}_{{G_{i}},k}^0(P)} $. If $u_x\notin P$, then  $u_x \in {{\mathscr P}_{{S_p^3},k}^t(P)}$ for some integer $t>0$. Note that in $S_p^3$, $u_x$ cannot be observed until one of its neighbors is. Since $u_x=a^3$, its $p-1$ neighbors in $S_p^3$ have labels in the form $aab$ for $b=0,\ldots ,p-1$, $b\not= a$. Therefore, $N_{S_p^3}[u_x]$ induces a $p$-clique in $S_p^3$. This means that when a neighbor forces $u_x$, it also forces all the vertices in the $p$-clique induced by $N_{S_p^3}[u_x]$. When this happens in $G_{i}$ instead of in $S_p^3$, $u_x$ has exactly one  unobserved neighbor ($x$), so $x \in {{\mathscr P}_{{S_p^3},k}^{t+1}(P)}$. 

Finally we prove $\gamma_{P,k}(G_{i})\geq \gamma_{P,k}({S_p^3})$ by showing that there exists a minimum $k$-power dominating set $Q$ of $G_{i}$ that is also a $k$-power dominating set of $S_p^3$. By Lemma \ref{highdegkPSD} there exists a minimum $k$-power dominating set $Q$ of $G_{i}$ that does not contain vertices of degree $1$, so $Q\subseteq  V({S_p^3})$. Therefore, in the $k$-power domination process starting with $Q$ in $G_{i}$, a vertex of degree $1$ in $G_{i}$ cannot be observed until its one neighbor in $S_p^3$ is. Then $Q$ is a $k$-power dominating set of $S_p^3$. 
We conclude that $\gamma_{P,k}(\widehat {sS_p^3})= \gamma_{P,k}(S_p^{3})$.  
\epf

It is known that if $n\geq 3$, $k\geq 1$, and $p\geq k+2$, then $\gamma_{P,k}(S_p^{n})=p^{n-2}(p-k-1)$ \cite{DK14} so we immediately obtain the following result.

\begin{cor}\label{sier}
Given integers  $n\ge 4$, $k\geq 1,$ and $p\geq k+2$, let  $s$ be a $(n-3)$-tuple of integers in $\{0,\ldots ,p-1\}$. Then $\gamma_{P,k}(\widehat {sS_p^3})=p(p-k-1)$.
\end{cor}

\begin{lem}\label{sierp}
Given integers $n\ge 4$, $k\geq 1$, and $p\geq k+2$, let $T$ denote the set of all  $(n-3)$-tuples of integers in $\{0,\ldots ,p-1\}$.  Then  $\gamma_{P,k}(S_p^n)= \sum _{t\in T}\gamma_{P,k}(\widehat {tS_p^3})$.
\end{lem}

\bpf
By Lemma \ref{sier} $\gamma_{P,k}(\widehat {tS_p^3})= \gamma_{P,k}(S_p^3)=p(p-k-1)$ for every $t\in T$. There are $p^{n-3}$ tuples in $T$, so $\sum _{t\in T}\gamma_{P,k}(\widehat {tS_p^3})= p^{n-3}p(p-k-1)=p^{n-2}(p-k-1)=\gamma_{P,k}(S_p^n)=p^{n-2}(p-k-1)$. 
\epf

The bound in Theorem  \ref{partitionpd} is $\gamma_{P,k}(S_p^n)\leq \sum _{t\in T}\gamma_{P,k}(\widehat {{tS_p^3}})$, so the following result is an immediate consequence of Lemma \ref{sierp}.

\begin{cor}
The bound in Theorem \ref{partitionpd} is tight.
\end{cor}

Next we present the equivalent results for $k$-forcing taking into consideration the following differences between $k$-power domination and $k$-forcing. The proofs are analogous and are omitted.

 \begin{itemize}
 \item[1.] For the lower bound, note that if $\deg_H(v_X)> k$, $\{v_X\}$ does not force $N_H(v_X)=N_G[X]\setminus X$. In that case, to obtain a $k$-forcing set of $H$ from a $k$-forcing set of $G$ it might be necessary to add at most $|N_G[X]\setminus X|-k$ vertices.
 \item [2.] For the upper bound, since there is no $k$-forcing equivalent to Lemma \ref{3}, it could happen that every minimum $k$-forcing set of $\widehat X$ contains a vertex $x\in N[X]\setminus X$ for which $\deg_{\widehat X}(x)=1$ but   $\deg_{G}(x)>k$. Thus, $x$ forces its neighbors in $\widehat X$ but not in $G$, and a $k$-forcing set of $\widehat X$ might not force $X$ in $G$.
 \end{itemize}

\begin{prop}\label{subgraphzf}
Let $G=(V,E)$ be a connected graph. Let $X\subseteq V(G)$. If there exists a minimum $k$-forcing set of $\widehat X$ that contains only vertices in $X$, then

\[Z_{k}(G/X)+Z_k(\widehat X)\geq Z_{k}(G) \geq \left\{ \begin{array}{ll}
      Z_{k}(G/X)-1&\mbox{ if } |N[X]\setminus X|\leq k, \\
       Z_{k}(G/X)- |N[X]\setminus X|+k &\mbox{ if }  |N[X]\setminus X|> k. 
     \end{array} \right.\]
     
\end{prop}

\begin{prop}\label{contractionubzf} Let $G=(V,E)$ be a connected graph.  Let $X\subseteq V$ such that  $\deg_G(x)\leq k$ for $x\in X$. If there exists a minimum $k$-forcing set of $\widehat X$ that contains only vertices in $X$, then
        
\[Z_{k}(G/X)+c(G[X])\geq Z_{k}(G) \geq \left\{ \begin{array}{ll}
      Z_{k}(G/X)-1&\mbox{ if } |N[X]\setminus X|\leq k,\\
       Z_{k}(G/X)- |N[X]\setminus X|+k &\mbox{ if }  |N[X]\setminus X|> k. 
     \end{array} \right.\]
 \end{prop}

\begin{cor}{\rm  \cite[Theorem 5.1]{O09}}
For any edge $e$ in a graph $G$, $Z(G)-1\leq Z(G/e)\leq Z(G)+1.$ 
\end{cor}

\begin{thm}\label{partitionkf}
Let $G=(V,E)$ be a  connected graph and let $P_1,\ldots ,P_r$ be a partition of $V$. If  $\widehat {P_i}$ has a minimum $k$-power dominating set in $P_i$ for every $i=1,\ldots ,r$, then

\[\Z _{k}(G)\leq \sum_{i=1}^r \Z_{k}(\widehat {P_i}).\]
\end{thm}
 
Theorem \ref{partitionpd} and Theorem \ref{partitionkf} provide upper bounds for the $k$-power domination and the $k$-forcing number of a graph in terms of  the $k$-power domination and the $k$-forcing number of  $\widehat {P_1},\ldots , \widehat {P_r}$, which can be computed in parallel. In particular, the importance of Theorems \ref{partitionpd} and \ref{partitionkf} resides in the fact that  $\widehat P_i$ might have properties that do not hold for $G$. For example, suppose $G$ is not a tree, but there is a linear algorithm to partition $V(G)$ into sets $P_1,\ldots ,P_r$ such that $\widehat{P_1},\ldots ,\widehat{P_r}$ are trees. Then using the linear algorithm for trees provided in \cite{DHLMR13}, $\gamma _{P,k}(\widehat X)$ can be computed in linear time. 
The exploration of possible uses of our results in algorithms to find $k$-power dominating or $k$-forcing sets a graph requires a detailed and careful analysis that is outside the scope of this paper. 


\section{$k$-power domination and $k$-forcing numbers}

By Observation \ref{basic},  $\gamma _{P,k}(G)\leq Z_k(G)\leq \gamma _{P,k}(G) (\Delta (G)+1). $ In this section we improve the upper bound in the previous inequality by generalizing a result by Benson et al.  \cite[Theorem 3.2]{REUF2015} for 1-power domination. An important concept in this work is that of {\it private neighborhood}, which we recall. Suppose $v\in S\subseteq V$.  A {\em $S$-private neighbor} of $v$ is a vertex $x \in N(v)$ such that $x \notin  N(S\setminus\{v\})$. Moreover, we say that $x$ is an {\em external} $S$-private neighbor if $x\notin S$. 

\begin{lem}\label{pvtnhbrs}{\rm \cite[Lemma 10]{CDMR12}} 
 In every connected graph $G$ with $\Delta (G) \geq k+2$ there exists a minimum $k$-power dominating set $S$ in which every vertex has at least $k+1$ $S$-private neighbors.  
 \end{lem}

We strengthen Lemma \ref{pvtnhbrs}  by extending it to external private neighbors.

\begin{lem}\label{epn}
 In every connected graph $G$ with $\Delta (G) \geq k+2$ there exists a minimum $k$-power dominating set $S$ in which every vertex has at least $k+1$ external $S$-private neighbors.  
 \end{lem}
 
 \bpf
By Lemma \ref{pvtnhbrs} there exists a minimum $k$-power dominating set $S$ in which every vertex has at least $k+1$ $S$-private neighbors. Suppose that there exists $u\in S$ that has at most $k$ external private neighbors.  We prove that $S'=S\setminus \{u\}$ is $k$-power dominating set, which contradicts the minimality of $S$. Since $u$ has at least $k+1$ neighbors and at most $k$ of them are outside $S$, there exists $y\in S$ such that $u$ and $y$ are neighbors. This implies that $u\in {{\mathscr P}_{G,k}}^0(S')$. Moreover, all non-external neighbors of $u$ are in $S$ and thus in ${{\mathscr P}_{G,k}}^0(S')$ so $u$ has at most $k$ unobserved neighbors.  Thus, $N_G[u]\subseteq {{\mathscr P}_{G,k}}^1(S')$, 
 so $S'$ is a $k$-power dominating set.
 \epf

\begin{lem}\label{Deanetal-k}
	 If $G$ is a connected graph with $\Delta  (G) \geq k+2$ and $S=\{u_1,\dots,u_t\}$ is a minimum $k$-power dominating set of $G$ in which every vertex has at least $k+1$ external $S$-private neighbors, then  \[\Z_k(G)\le \sum_{i=1}^t(\deg u_i+1-k).\]
\end{lem}

\bpf By hypothesis,  for each $i=1,\ldots ,t$ there exists a set $\{ x_1^{(i)},\dots,x_k^{(i)} \}$ of external $S$-private neighbors of $u_i$. We prove that
$B:= \bigcup_{i=1}^t \left(N[u_i]\setminus\{x_1^{(i)},\dots,x_k^{(i)}\}\right)$ is a $k$-forcing set of $G$. Since $x_1^{(i)},\dots,x_k^{(i)}$ are external $S$-private neighbors of $u_i$, then $\{ x_1^{(i)},\dots,x_k^{(i)} \}\cap S=\emptyset$, which implies $u_i\in B$, for every $i=1,\ldots ,t$. In the first step of the $k$-forcing process each vertex $u_i$ forces $x_1^{(i)},\dots,x_k^{(i)}$ so $B$ is a $k$-forcing set of $N[S]$ in $G$.  Since $S$ is a $k$-power dominating set of $G$, by Observation \ref{basic} $N[S]$ is a $k$-forcing set of $G$. Then by Observation \ref{transitive},  $B$ is a $k$-forcing set of $G$ so,  
$\Z_k(G)\le |B|\le \sum_{i=1}^t \left| N[u_i]\setminus\{x_1^{(i)},\dots,x_k^{(i)}\}\right|\le \sum_{i=1}^t(\deg u_i+1-k).$ 
\epf

\begin{thm}\label{REUFbound}
  In every connected graph $G$ with $\Delta (G)\geq k+2$,   \[ \Z_k(G)\leq  \pdk k(G) ({\Delta(G)+1-k}), \mbox{ or equivalently, }   \lc\frac{\Z_k(G)}{\Delta(G)+1-k}\rc\le \pdk k(G)\] and this lower bound for  $\gamma _{P,k}(G)$ is tight.
\end{thm}

\bpf  
     By Lemma \ref{epn} there exists a minimum $k$-power dominating set $S=\{u_1,\dots,u_{\pdk k(G)}\}$ of $G$ in which each vertex has at least $k+1$ external $S$-private neighbors.  By Lemma \ref{Deanetal-k}, $\Z_k(G)\leq \sum_{i=1}^{\pdk k(G)}(\deg u_i+1-k)\le \pdk k(G)(\Delta(G) +1-k)$, and as a consequence, $ \lc\frac{\kf (G)}{\Delta(G)+1-k}\rc\leq \kpd (G).$

To prove that the bound is tight, let $r\ge 2$ and  $p>3r+k-3$. Construct the graph $G_{p,r}$ by adding $p$ pendent vertices to each vertex $v_i$, $i=1,\ldots ,r$ to each vertex of a path of order $r$. Then $\Delta (G_{p,r})=p+2$ and since $p\geq k+1$,  $ \pdk k(G_{p,r})=r$, $\Z_k(G_{p,r})=r(p-k)$, and $\lc\frac{r(p-k)}{p+3-k}\rc=r$.
\epf

Next we apply Theorem \ref{REUFbound} to obtain lower bounds for the $k$-forcing number of graphs from upper bounds for the $k$-power domination number of an arbitrary graph presented in  \cite {CDMR12}  and improved in \cite {DHLMR13} for $(k+2)$-regular graphs.

\begin{thm}\label{uppd}{\rm  \cite [Theorem 11] {CDMR12}} Let $G$ be a connected graph with  $|G|\geq k+2$. Then $\gamma _{P,k}(G) \leq {|G|\over {k + 2}}$.
\end{thm}

\begin{cor} \label{lbzf} In a connected graph $G$ with  $\Delta (G) \geq k+2$, \[\Z _k(G) \leq \lf {{|G|\over {k + 2}}(\Delta (G) +1-k)}\rf \] and this bound is tight.
\end{cor}

\bpf Since $\Delta (G) \geq k+2$ implies $|G|\geq k+2$ we apply  Theorem \ref{uppd} and obtain   $\gamma _{P,k}(G) \leq {|G|\over {k + 2}}$. By Theorem \ref{REUFbound} we know $\Z _k(G)\leq \gamma _{P,k} (G)(\Delta (G)+1-k)$ and combining both inequalities we conclude $\Z _k(G) \leq \lf{ {{|G|}\over {k + 2}}(\Delta (G)+1-k)}\rf$.

To show this bound is tight, observe that $\Z_k(K_{k+3})=3$, and the upper bound in this case is $\lf {{k+3}\over {k+2}}(k+2+1-k)\rf=3$ for $k\geq 2$.
\epf

\begin{thm}{\rm \cite [Theorem 2.1] {CDMR12}} Let $G$ be a connected $(k +2)$-regular graph. If $G \not= K_{k+2,k+2}$, then $\gamma _{P,k}(G) \leq {|G|\over {k + 3}}$.
\end{thm}

\begin{cor} \label{lbzfr}  Let $G$ be a connected $(k +2)$-regular graph.  If $G \not= K_{k+2,k+2}$, then $\Z _k(G) \leq {{3|G|}\over {k + 3}}$.
\end{cor}

\bpf
Since $G$ is $(k+2)$-regular,  $\Delta (G)= k+2$ so we apply Theorem \ref{REUFbound} and obtain  $\Z _k(G)\leq {|G|\over {k + 3}}(k+2+1-k)= {{3|G|}\over {k + 3}}$.  To see that the bound is best possible it suffices to consider 
$ K_{k+3}$ which is $(k+2)$-regular and $\Z _k(K_{k+3})=3= \frac{3(k+3)} {k + 3}$. 
\epf


\end{document}